

\baselineskip=14pt
\parskip=10pt

\font\eightrm=cmr8 

\magnification=\magstephalf

\def\1{{\overline{1}}}
\def\2{{\overline{2}}}
\parindent=0pt
\overfullrule=0in

\def\frac#1#2{{#1 \over #2}}
\centerline
{\bf Some Remarks on a recent article by J. -P. Allouche }
\bigskip
\centerline
{\it Shalosh B. EKHAD and Doron ZEILBERGER}
\bigskip

In a recent article [A], the author tries very hard to ``explain" a puzzling remark made by  Otto G. Ruehr in his
solution [R] to American Mathematical Monthly elementary problem {\bf E2765} (proposed in 1979 by Naoki Kimura).

Frankly, we don't see the point of going to such great lengths to prove that routinely-provable identity $A$
implies routinely-provable identity $B$. It is much faster to  prove them {\bf both} from scratch.
As we will soon see, Ruehr's remark was most probably a {\it non-sequitur}, and it is very unlikely that
it followed [A]'s {\it exegesis}.

The original problem was to prove that for every continuous function $f$ the following identity holds
$$
\int_{-\frac{1}{2}}^{\frac{3}{2}} f(3x^2-2x^3) \, dx \, = \,
2 \, \int_{0}^{1} f(3x^2-2x^3) \, dx  \quad .
$$

After giving his change-of-variable proof, Ruehr makes the following remark.

{\it ``An interesting alternative approach would be to use the Weierstrass Approximation Theorem to reduce that problem
to that of a polynomial $f$, and in turn, by linearity to that of establishing the given equation for $f(z)=z^n$.''}

In other words one had to prove, that for every non-negative integer $n$,
$$
\int_{-\frac{1}{2}}^{\frac{3}{2}} (3x^2-2x^3)^n \, dx \, = \,
\int_{0}^{1} 2 \, (3x^2-2x^3)^n \, dx  \quad .
\eqno(Ruehr)
$$

While this was non-trivial back in 1980, it is routinely provable today, using the
{\it Almkvist-Zeilberger} algorithm [AZ], implemented in the Maple package {\tt EKHAD.txt} [Z].
Just type:

{\tt AZdI((3*x**2-2*x**3)**n,x,n,N,-1/2,3/2)[1];} \quad {\tt AZdI(2*(3*x**2-2*x**3)**n,x,n,N,0,1)[1];}

for the left and right sides, respectively, and in a split second you would get the {\bf same} output:

{\tt [9 (n + 1) (2 n + 1) - 2 (3 n + 4) (3 n + 2) N, 2]} \quad ,

which means that {\bf both} the left side and the right side satisfy the first-order inhomogeneous linear recurrence equation
with polynomial coefficients
$$
9(n+1)(2n+1)f(n) \, - \, 2(3n+4)(3n+2)f(n+1) \, = \, 2 \quad .
$$
Since $L(0)=R(0)=2$, this immediately implies $(Ruehr)$, without any need for a clever change of variable.

If you take Ruehr's suggestion literally, then you would get a certain binomial coefficients identity that is easily
provable by the Zeilberger algorithm, supplying yet another routine proof. 
This summation identity, that the reader can easily derive by using
the  binomial theorem on $(3-2x)^n$ and integrating term-by-term (as suggested in [R])
has {\bf nothing whatsoever to do} with the following two binomial coefficients identities
that Ruehr claims are `equivalent'. We would never know what Ruehr had in mind (he probably
got mixed up with another problem), and frankly we don't really care.
At any rate, we are sure that it is not via Allouche's extremely circuitous route.

These `equivalent' (per Ruehr) identities are 
$$
\sum_{j=0}^{2n} (-4)^j {{3n+1} \choose {n+j+1}} \, = \,\sum_{j=0}^{n} 2^j {{3n+1} \choose {n-j}} \quad ,
$$
$$
\sum_{j=0}^{2n} (-3)^j {{3n-j} \choose {n}} \, = \,\sum_{j=0}^{n} 3^j {{3n-j} \choose {2n}} \quad .
$$

As already pointed out in [MTWZ] these two new identities (that, in spite of [A], are  completely unrelated to the original
integral identity $(Ruehr)$)
are also {\bf routinely provable} with Maple. In fact it is also strange that they are listed as two different identities.
All four sums happen to be identical.

Go into Maple, and type

{\tt Z:=SumTools[Hypergeometric][ZeilbergerRecurrence];} \quad \quad \quad \quad \quad ,

then type

{\tt  Z(3**j*binomial(3*n-j,2*n),n,j,f,0..n);} \quad {\tt  Z((-3)**j*binomial(3*n-j,n),n,j,f,0..2*n);}

{\tt  Z(2**j*binomial(3*n+1,n-j),n,j,f,0..n);} {\tt  Z((-4)**j*binomial(3*n+1,n+j+1),n,j,f,0..2*n);}

and find out, in one nano-second, that all four sums satisfy the {\bf same} linear recurrence, namely
$$
-27\,f \left( n \right) +4\,f \left( n+1 \right) \, = \,
-3\, \frac{(3n+1)!}{(2n+1)!(n+1)!}
$$

As also pointed out in [MTWZ] these four sums are all equal to OEIS sequence  {\tt A6256} [S] whose definition is yet another binomial coefficient
sum
$$
\,\sum_{k=0}^{n} {{3k} \choose {k}}  {{3n-3k} \choose {n-k}} 
\quad .
$$
Typing 

{\tt Z(binomial(3*k,k)*binomial(3*n-3*k,n-k),n,k,f,0..n); }

gives that it satisfies a second-order {\bf homogeneous} linear recurrence equation with polynomial coefficients
$$
-81\, \left( 3\,n+2 \right)  \left( 3\,n+4 \right) f \left( n \right) + \left( 216\,{n}^{2}+594\,n+420 \right) f \left( n+1 \right) -8\, \left( 2\,n+3
 \right)  \left( n+2 \right) f \left( n+2 \right) \, = \, 0 \quad ,
$$
that is routinely equivalent to the above first-order inhomogeneous recurrence for the Ruehr sums.

{\bf Conclusion}

In the conclusion to [A], the author states

{\it ``The literature about sums involving binomial coefficients is huge.''} 

He should have added, 

{\it ``... and mostly obsolete\footnote{$^1$}
{\eightrm There is still room for {\bf elegant} and insightful combintorial  proofs of such identities, of course, see [MTWZ]
for two such gorgeous proofs.}
 (thanks to the Wilf-Zeilberger [PWZ] algorithmic proof theory, implemented in Maple, Mathematica and other systems)''}.

Both the integral formula with $f(z)=z^n$ and the two binomial coefficients identities, are nowadays routinely provable,
and it has very little {\it mathematical} interest to do {\it exegesis} of what Ruehr had in mind.
Of course, it may be of some {\it psychological} or {\it literary} interest. After all literary scholars often
go to great lengths to try and understand `what the great poet had in mind', and most often get it wrong.
Often the poet had nothing in mind, and if she did, it was something very mundane.

{\bf Parody}

Let's use an analogy. Suppose that in the 
American Mathematical Monthly analog of Egypt in 5000BC there was a papyrus that gives an elegant proof of the identity
$$
123 \cdot 321 \, = \,  39483 \quad ,
\eqno(1)
$$
and then it comments that it is equivalent to
$$
111 \cdot 449 \, = \,    49839  \quad ,
\eqno(2)
$$
but that the latter identity is no easier to prove than the former one. In a logical sense, that ancient savant would have been
right, since all correct statements are logically equivalent. 
But some mathematical historian could have been tempted  to `explain' how Professor Ahmes may have reasoned.
Since $111= (123-12)$ and  $449=321 +128$ (proofs left to the reader), we have
$$
111 \cdot 449 \, = \, (123-12)(321 +128) \, = \,
123 \cdot 321 + 123 \cdot 128 - 12 \cdot 321 - 12 \cdot 128
$$
$$
=123 \cdot 321 + 15744 - 3582-1536 \quad.
$$
Using Eq. $(1)$, the first product equals  $39483$, and  it follows that indeed
$$
111 \cdot 449 \,= \, 39483 + 15744 - 3852-1536 \, = \, 49839 \quad  QED.
$$
With all due respect, it is much easier to prove identities $(1)$ and $(2)$ separately, than 
to only prove $(1)$ and then to deduce $(2)$ from it.

{\bf References}

[A] J. -P. Allouche, {\it Two binomial identities of Ruehr Revisited}, Amer. Mathematical Monthly, {\bf 126}(3) [March 2019], 217-225.

[AZ] Gert Almkvist and Doron Zeilberger, {\it The method of differentiating under the integral sign}
J. Symbolic Computation {\bf 10} (1990), 571-591. Available from \hfill\break
{\tt http://sites.math.rutgers.edu/\~{}zeilberg/mamarim/mamarimhtml/duis.html}   \quad .

[MTWZ] Sean Meehan, Akalu Tefera, Michael Weselcouch, and Aklilu Zeleke, {\it Proof os Ruehr's identities},
INTEGERS {\bf 14} (2014), \#A10 (6 pages) \quad .

[PWZ] Marko Petkovsek, Herbert S. Wilf, and Doron Zeilberger, {\it ``A=B''}, A.K. Peters, 1996. Freely available from \hfill\break
{\tt https://www.math.upenn.edu/\~{}wilf/AeqB.html} \quad .

[R] Otto G. Ruehr, {\it  Solution to Problem E2765 (proposed by Naoki Kimura)},
Amer. Mathematical Monthly, {\bf 87}(4) [April 1980], 307-308.

[S] Neil Sloane, {\it OEIS sequence A006256} {\tt http://oeis.org/A006256}.

[Z] Doron Zeilberger, {\bf EKHAD.txt}, a Maple package freely available from \hfill\break
{\tt http://sites.math.rutgers.edu/\~{}zeilberg/tokhniot/EKHAD.txt}  \quad .

\bigskip
\hrule
\bigskip
Shalosh B. Ekhad, c/o D. Zeilberger, Department of Mathematics, Rutgers University (New Brunswick), Hill Center-Busch Campus, 110 Frelinghuysen
Rd., Piscataway, NJ 08854-8019, USA. \hfill\break
Email: {\tt ShaloshBEkhad at gmail dot com}   \quad .
\bigskip
Doron Zeilberger, Department of Mathematics, Rutgers University (New Brunswick), Hill Center-Busch Campus, 110 Frelinghuysen
Rd., Piscataway, NJ 08854-8019, USA. \hfill\break
Email: {\tt DoronZeil at gmail  dot com}   \quad .
\bigskip
\hrule
\bigskip
Exclusively published in the Personal Journal of Shalosh B.  Ekhad and Doron Zeilberger and arxiv.org \quad .
\bigskip
Written: Purim 5779 (alias March 21, 2019).
\end